\documentclass[number,citesort,seceqn,dvips]{arxbj}
\usepackage{cursive}  
\usepackage{graphicx}

\aid{0}
\volume{17}
\issue{1}
\pubyear{2011}
\firstpage{484}
\lastpage{506}
\doi{10.3150/10-BEJ281}

\makeatletter

\newtheorem{theorem}{Theorem}[section]
\newtheorem{lemma}[theorem]{Lemma}

\newremark{example}[theorem]{Example}
\newtheorem{proposition}[theorem]{Proposition}
\newproclaim{definition}[theorem]{Definition}
\newtheorem{corollary}[theorem]{Corollary}

\def\eqref#1{(\ref{#1})}

\makeatother

\begin{document}
\begin{frontmatter}

\title{Fractional L\'evy-driven Ornstein--Uhlenbeck processes and
stochastic differential equations}
\runtitle{Fractional L\'evy-driven Ornstein--Uhlenbeck processes and SDEs}

\begin{aug}
\author[a]{\fnms{Holger} \snm{Fink}\thanksref{a}\corref{}\ead[label=e1]{fink@ma.tum.de}\ead[label=u1,url]{http://www-m4.tum.de/pers/fink}} \and
\author[b]{\fnms{Claudia} \snm{Kl\"uppelberg}\thanksref{b}\ead[label=e2]{cklu@ma.tum.de}\ead[label=u2,url]{http://www-m4.ma.tum.de}}
\runauthor{H. Fink and C. Kl\"uppelberg}
\address[a]{Center for Mathematical Sciences, Technische Universit\"at
M\"unchen, D-85748 Garching,
Germany. \printead{e1,u1}}
\address[b]{Center for Mathematical Sciences and Institute for
Advanced Study, Technische
Universit\"at M\"unchen, D-85748 Garching, Germany.
\printead{e2,u2}}
\end{aug}

\received{\smonth{8} \syear{2009}}
\revised{\smonth{4} \syear{2010}}

%
\begin{abstract}
Using Riemann--Stieltjes methods for integrators of bounded
$p$-variation we define a pathwise integral driven by a
fractional L\'evy process (FLP). To explicitly solve general fractional
stochastic differential equations
(SDEs) we introduce an Ornstein--Uhlenbeck model by a stochastic
integral representation, where the driving
stochastic process is an FLP. To achieve the convergence of improper
integrals, the long-time behavior of FLPs
is derived. This is sufficient to define the fractional L\'evy--Ornstein--Uhlenbeck process (FLOUP) pathwise as an
improper Riemann--Stieltjes integral. We show further that the FLOUP is
the unique stationary solution of the
corresponding Langevin equation. Furthermore, we calculate the
autocovariance function and prove that
its increments exhibit long-range dependence.
Exploiting the Langevin equation, we consider SDEs driven by FLPs of
bounded $p$-variation for $p<2$ and
construct solutions using the corresponding FLOUP. Finally, we consider
examples of such SDEs, including
various state space transforms of the FLOUP and also fractional L\'evy-driven Cox--Ingersoll--Ross (CIR) models.
\end{abstract}

%

\end{frontmatter}

\section{Introduction}\label{s1}

In this paper we consider (stationary) solutions to SDEs of the form
%
\begin{eqnarray}\label{qq11}
\mathrm{d}X_t=\mu(X_t)\,\mathrm{d}t+\sigma(X_t)\,\mathrm{d}L^d_t,\qquad t\in\mathbb{R},
\end{eqnarray}
where $L^d$ is a fractional L\'evy process (FLP) of bounded
$p$-variation for $p<2$ and $\mu$ and $\sigma$ are appropriate
coefficient functions. Applying pathwise Riemann--Stieltjes integration
for functions of bounded $p$-variation, we solve such equations by
constructing explicit solutions. The basic model will be a fractional
L\'evy--Ornstein--Uhlenbeck process (FLOUP) introduced by the
stochastic integral representation
\begin{eqnarray*}
\mathcal{L}^{d,\lambda}_t=\int_{-\infty}^t\mathrm{e}^{-\lambda(t-u)}\,\mathrm{d}L^d_u
,\qquad t\in\mathbb{R}.
\end{eqnarray*}
We further show that this is the unique stationary pathwise solution of
the corresponding Langevin equation
%
\begin{eqnarray}\label{lange}
\mathrm{d}\mathcal{L}^{d,\lambda}_t=-\lambda\mathcal{L}^{d,\lambda}_t \,\mathrm{d}t +
\mathrm{d}L^d_t ,\qquad t\in\mathbb{R}.
\end{eqnarray}
Using this relation we will consider SDEs of the form $(\ref{qq11})$
and impose assumptions on the coefficient functions $\mu$ and $\sigma
$, under which solutions can be constructed by monotone transformation
of $\mathcal{L}^{d,\lambda}$.

Although our paper is purely theoretical, we are aiming at applications
where positive solutions of \eqref{qq11} are of interest.
An approach, developed in \cite{KL} for
SDEs driven by FBM, can be modified to SDEs driven by FLPs.
On the other hand, a squared FLOUP is positive and a solution to the SDE
\begin{eqnarray*}
\mathrm{d}X_t=-2\lambda X_t \,\mathrm{d}t + 2\sqrt{|X_t|}\,\mathrm{d}L^d_t ,\qquad t\in\mathbb{R}.
\end{eqnarray*}
We will discuss various examples with different state spaces and
different $\mu$ and $\sigma$.
We will also present some properties of the respective solutions, also
concerning the stationary distribution.

Our paper is organized as follows.
Section~\ref{sec2} considers FLPs and pathwise integration. Section~\ref{sec3} introduces
the FLOUP as a pathwise improper Riemann--Stieltjes integral and shows
that it is the unique stationary pathwise solution of the corresponding
Langevin equation. Moreover, we calculate its autocovariance function
and show that the increments of an FLOUP exhibit long-range dependence.
Section~\ref{sec4} mainly extends Buchmann and Kl\"uppelberg \cite{KL} from
fractional Brownian motion to FLPs and states structural conditions for
the coefficient functions $\mu$ and $\sigma$, which guarantee an
existence (and uniqueness) theorem. Section~\ref{sec5} provides examples and
simulations. The \hyperref[app]{Appendix} reviews the Riemann--Stieltjes analysis via
$p$-variation.

The following notation will be used throughout. We always assume a
complete probability space $(\Omega,\mathcal{F},P)$. We denote the
$\mathcal{F}$-measurable real functions by $L^0(\Omega)$, the Hilbert
space of square integrable random variables by $L^2(\Omega)$, the
vector space of continuous real functions on $A\subseteq\mathbb{R}$
by $\mathcal{C}^0(A)$ and by $\|\cdot\|_{\sup}^{A}$ the supremum norm.
Furthermore, $\operatorname{Lip}(A)$ and $C^1(A)$ are the spaces of real
functions on $A$, which are Lipschitz continuous on compacts and
continuously differentiable, respectively.
The spaces of integrable and square integrable real functions are
denoted by $L^1(\mathbb{R})$ and $L^2(\mathbb{R})$, respectively.
When speaking of a two-sided L\'evy process we mean the following: given
two independent copies of the same L\'evy process, $L^1$ and $L^2$, we take
%
\begin{eqnarray}\label{qwas1}
L_t:=L^1_t1_{\{t\geq0\}}+L^2_{-t-}1_{\{t<0\}},\qquad t\in\mathbb{R}.
\end{eqnarray}
The Dirac measure in $1$ we denote by $\delta_1$.
Finally, for $-\infty\leq b\leq a\leq\infty$ we set $[a,b]:=[b,a]$.

Integrals throughout this paper are considered in the
Riemann--Stieltjes sense, if not stated otherwise.

\section{Fractional L\'evy processes and pathwise integration}\label{sec2}

Fractional L\'evy processes (FLPs) were introduced as a natural
generalization of the integral representation of fractional Brownian
motion (FBM). We shortly review the main properties of FLPs, see
\cite{tina}, Section~\ref{sec3}, for details and more background.
For notational convenience we work with the fractional integration
parameter $d\in(-\frac{1}{2},\frac{1}{2})$ instead of the Hurst
index $H$, where $d=H-\frac{1}{2}$.
Because we are only interested in long memory models, we restrict
ourselves to $d\in(0,\frac{1}{2})$.
Furthermore, we only consider FLPs with existing second moments.
Analogously to \cite{mandel} for FBM we choose
(like Marquardt~\cite{tina}) the following definition.

\begin{definition}\label{def2.1}
Let $L=(L_t)_{t\in\mathbb{R}}$ be a zero-mean two-sided L\'evy
process with
$E[L(1)^2]<\infty$ and without a Brownian component.
For $d\in(0,\frac{1}{2})$ we define
%
\begin{equation}\label{FLP2}
L^d_t:=\frac{1}{\Gamma(d+1)}\int_{-\infty}^\infty
[(t-s)_+^d-(-s)_+^d]L(\mathrm{d}s),\qquad t\in\mathbb{R}.
\end{equation}
We call $L^d=(L^d_t)_{t\in\mathbb{R}}$ a
\emph{fractional L\'evy process (FLP)} and $L$ the \emph{driving L\'
evy process} of $L^d$.
\end{definition}

The integrals above exist in the $L^2(\Omega)$-sense; see
\cite{tina}, Theorem~3.5, for details.

Recall that, by the L\'evy--It\^o decomposition, every L\'evy process
can be represented as the sum of a Brownian motion and an independent
jump process.
The Brownian motion gives rise to an FBM, which has been studied
extensively; see, for instance, \cite{samota}
for general background or \cite{KL} in the
context of the present paper.

The next result ensures that there is, in fact, a modification of
$(\ref{FLP2})$ that equals a pathwise improper Riemann integral and
gives first properties.

\begin{proposition}[(\cite{tina}, Theorems 3.7, 4.1 and
4.4)]\label{pr3.2}
Let $L^d$ be an FLP with $d\in(0,\frac{1}{2})$.
Then the following assertions hold:
\begin{enumerate}[(iii)]
\item[(i)] $L^d$ has a modification that equals the improper Riemann integral
%
\begin{equation}\label{FLP6}
\frac{1}{\Gamma(d)}\int_\mathbb
{R}[(t-s)_+^{d-1}-(-s)_+^{d-1}]L(s)\,\mathrm{d}s,\qquad t\in\mathbb{R}.
\end{equation}
Furthermore, $(\ref{FLP6})$ is continuous in $t$.
\item[(ii)] For $s,t\in\mathbb{R}$ we have
%
\begin{equation}\label{FLP7}
\operatorname{Cov}(L^d_t,L^d_s)=\frac{E[L(1)]^2}{2\Gamma(2d+2)\sin(\curpi
(d+{1/2}))}[|t|^{2d+1}+|s|^{2d+1}
-|t-s|^{2d+1}] .
\end{equation}
\item[(iii)] $L^d$ has stationary increments and is symmetric, i.e.,
$(L^d_{-t})_{t\in\mathbb{R}}\stackrel{d}{=}(-L^d_t)_{t\in\mathbb{R}}$.
\end{enumerate}
\end{proposition}

From now on, we always work with the modification of Proposition~\ref
{pr3.2}(i).

Next we define integration with respect to FLPs.
As has been shown in \cite{tina}, Theorem~4.10, FLPs may not
be semimartingales, and integration in the $L^2(\Omega)$-sense has
been developed in~\cite{tina}, Section~5.
Theorem 4.4 in~\cite{tina} also shows that FLPs are only H\"
older continuous up to the fractional integration parameter $d$ and not
to the Hurst index $H$ as in the case for FBM. Therefore, pathwise
Riemann--Stieltjes integration by H\"older continuity does not make
sense for SDEs.
On the other hand, using an approach like Young~\cite{Young} based on
$p$-variation of the sample paths, integration in a pathwise
Riemann--Stieltjes sense can be defined; for details see the \hyperref[app]{Appendix}.
This means we have a chain rule and a density formula provided the
integrator is of bounded $p$-variation for $p\in[1,2)$.

We recall the definition of $p$-variation over a compact interval
$[a,b]\subset\mathbb{R}$.
Let $f\dvtx [a,b]\mapsto\mathbb{R}$.
We define for $0<p<\infty$ the $p$-variation of $f$ as
%
\begin{equation}\label{rev:1}
v_p(f,[a,b]):=\sup_\kappa\sum_{i=1}^n|f(x_i)-f(x_{i-1})|^p,
\end{equation}
where the supremum is taken over all subdivisions $\kappa$ of $[a,b]$.
If $v_p(f,[a,b])<\infty$, then we say that $f$ is of \emph{bounded
$p$-variation} on $[a,b]$. We will further call an FLP $L^d$ of bounded
$p$-variation if it is a.s. of bounded $p$-variation on compacts.

Let $L^d$ be an FLP of bounded $p$-variation, $d\in(0,\frac{1}{2})$
and $p\in[1,2)$. For $A\subset\mathbb{R}$ we define
%
\begin{eqnarray}\label{wcon}
\mathfrak{W}^{\mathrm{con}}_p(A):=\{f\in\mathcal{C}^0(A)\dvtx
v_p(f,[s,t])<\infty \,\forall  [s,t]\subseteq\mathbb{R}\}.
\end{eqnarray}
Then we define for every stochastic process with sample paths $H\in
\mathfrak{W}^{\mathrm{con}}_q(\mathbb{R})$ a.s. and for $p,q>0$ with
$p^{-1}+q^{-1}>1$ the integral
%
\begin{equation}\label{IF2}
\int_a^bH_s\,\mathrm{d}L^d_s,\qquad-\infty\leq a\leq b\leq\infty,
\end{equation}
pathwise in the Riemann--Stieltjes sense.

As stated in the \hyperref[app]{Appendix} the integral in (\ref{IF2}) always
exists on finite intervals $[a,b]$. We consider also improper
integrals, where $a=\infty$ or $b=-\infty$.
The existence of the tail integral has then to be justified.

For example, FLPs, where the driving L\'evy process is of finite
activity, are of bounded $p$-variation for all $p\geq1;$ cf.
Theorem~2.25 of \cite{tina}.

\section{Fractional L\'evy--Ornstein--Uhlenbeck processes}\label{sec3}

We introduce fractional L\'evy--Ornstein--Uhlenbeck processes (FLOUPs)
as improper Riemann--Stieltjes integrals and prove that they are
stationary solutions of the Langevin equation \eqref{lange}.
To show the existence of the improper Riemann--Stieltjes integral, we
first need some knowledge about the long-time behaviour of FLPs. A
similar result considering $t\to\infty$ has been proven by Muneya
Matsui (personal communication).

\begin{theorem}\label{th3.1}
Let $L^d$ be an FLP, $d\in(0,\frac{1}{2})$.
Then for all $\alpha>d+\frac{1}{2}$ we have
%
\begin{equation}\label{FLP24}
\lim_{t\to-\infty}\frac{|L^d_t|}{|t|^\alpha}=0 \qquad\mbox{a.s.}
\end{equation}
\end{theorem}

\begin{pf}
Without loss of generality we can assume that $t<0$. By the law of the
iterated logarithm (LIL) for L\'evy processes (cf. \cite{sato},
Proposition 48.9) we find a random variable $T$ and a constant $M>0$
such that a.s. for all $s<T$\vspace*{2pt}
%
\begin{equation}\label{FLP16b}
|L(s)|\leq M(2|s|\log\log|s|)^{{1/2}}.\vspace*{2pt}
\end{equation}
We can always make $T$ smaller and so we choose $T<-\mathrm{e}$.
For any such path we can assume that $t<T$ and estimate\vspace*{2pt}
\begin{eqnarray}\label{fu1}
\frac{1}{|t|^\alpha}|L^d_t|&=&\frac{1}{|t|^\alpha}
\frac{1}{\Gamma(d)}\bigg|\int_{-\infty}^\infty
[(t-s)_+^{d-1}-(-s)_+^{d-1}]L(s)\,\mathrm{d}s
\bigg|\nonumber
\\[-9pt]\\[-9pt]
&\leq&
\frac{1}{\Gamma(d)}\frac{1}{|t|^\alpha}\biggl(\int_{-\infty
}^t[(-s)^{d-1}-(t-s)^{d-1}]|L(s)|\,\mathrm{d}s+
\int_{t}^0(-s)^{d-1}|L(s)|\,\mathrm{d}s\biggr) .\qquad\nonumber\vspace*{2pt}
\end{eqnarray}
Therefore, it suffices to show that a.s.\vspace*{2pt}
%
\begin{eqnarray}\label{I1}
\lim_{t\rightarrow-\infty}\frac{1}{|t|^\alpha}\int_{-\infty
}^t[(-s)^{d-1}-(t-s)^{d-1}]|L(s)|\,\mathrm{d}s=0\vspace*{2pt}
\end{eqnarray}
and\vspace*{2pt}
\begin{eqnarray}
\label{I2}
\lim_{t\rightarrow-\infty}\frac{1}{|t|^\alpha}\int
_{t}^0(-s)^{d-1}|L(s)|\,\mathrm{d}s=0.\vspace*{2pt}
\end{eqnarray}
We start with \eqref{I1}. Using the LIL we get an upper bound of (\ref
{fu1}) as follows:\vspace*{2pt}
\begin{eqnarray}\label{FLP24c}
&&\frac{1}{|t|^\alpha}\int_{-\infty
}^t[(-s)^{d-1}-(t-s)^{d-1}]|L(s)|\,\mathrm{d}s\nonumber
\\[2pt]
&&\quad\leq\frac{M}{|t|^\alpha}\int_{-\infty
}^{-|t|}[(-s)^{d-1}-(-|t|-s)^{d-1}](2|s|\log\log|s|)^{
{1/2}}\,\mathrm{d}s\nonumber
\\[-9pt]\\[-9pt]
&&\quad=\frac{M|t|}{\mathrm{e}|t|^\alpha}\int_{-\infty
}^{-\mathrm{e}}[(-\mathrm{e}^{-1}|t|u)^{d-1}-(-|t|-\mathrm{e}^{-1}|t|u)^{d-1}]\nonumber
\\[2pt]
&&{}\hspace*{42pt}\qquad\times(2\mathrm{e}^{-1}|t||u|\log
\log(\mathrm{e}^{-1}|t||u|))^{{1/2}}\,\mathrm{d}u,\nonumber\vspace*{2pt}
\end{eqnarray}
where we have used in the last line the change of variable $\mathrm{e}^{-1}|t|u=s$.
Now note that for large $|t|$ and $|u|\geq \mathrm{e}$\vspace*{2pt}
\begin{eqnarray*}\label{FLP24b}
|t||u|\log\log(\mathrm{e}^{-1}|t||u|)&=&|t||u|\log\bigl(\log(\mathrm{e}^{-1}|t|)+\log |u|\bigr)\nonumber
\\[2pt]
&\leq&|t||u|\log\log|t|+|t||u|\log(1+\log|u|) .\nonumber\vspace*{2pt}
\end{eqnarray*}

\noindent
Combining $(\ref{FLP24b})$ with $|a+b|^{{1/2}}\leq|a|^{
{1/2}}+|b|^{{1/2}}$ for $a,b\in\mathbb{R}$ we get an upper
bound for $(\ref{FLP24c})$ by\vspace*{2pt}
\begin{eqnarray}\label{FLP24w}
&&\frac{M(2\mathrm{e}^{-1}|t|\log\log|t|)^{{1/2}}}{\mathrm{e}|t|^{\alpha
-d}}\int_{-\infty
}^{-\mathrm{e}}[(-\mathrm{e}^{-1}u)^{d-1}-(-1-\mathrm{e}^{-1}u)^{d-1}]|u|^{{1/2}}\,\mathrm{d}u\nonumber
\\
&&\qquad{}+\frac{M(2\mathrm{e}^{-1}|t|)^{{1/2}}}{\mathrm{e}|t|^{\alpha-d}}\int_{-\infty
}^{-\mathrm{e}}[(-\mathrm{e}^{-1}u)^{d-1}-(-1-\mathrm{e}^{-1}
u)^{d-1}]\bigl(|u|\log(1+\log|u|)\bigr)^{{1/2}}\,\mathrm{d}u\qquad\nonumber
\\[-8pt]\\[-8pt]
&&\quad =\frac{M(2\mathrm{e}^{-1}\log\log|t|)^{{1/2}}}{\mathrm{e}|t|^{\alpha
-(d+{1/2})}}\int_{\mathrm{e}}^{\infty}[(\mathrm{e}^{-1}u)^{d-1}-(-1+\mathrm{e}^{-1}
u)^{d-1}]u^{{1/2}}\,\mathrm{d}u\nonumber
\\
&&\qquad{}+\frac{M(2\mathrm{e}^{-1})^{{1/2}}}{\mathrm{e}|t|^{\alpha-(d+
{1/2})}}\int_{\mathrm{e}}^{\infty}[(\mathrm{e}^{-1}u)
^{d-1}-(-1+\mathrm{e}^{-1}u)^{d-1}]\bigl(u\log(1+\log u)\bigr)^{{1/2}}\,\mathrm{d}u .\nonumber
\end{eqnarray}
By a binomial expansion we get
$(\mathrm{e}^{-1}u-1)^{d-1}=(\mathrm{e}^{-1}u)^{d-1}-(d-1)(\mathrm{e}^{-1}u)^{d-2}+\mathrm{O}(u^{d-3})$
and, therefore (writing $a(u)\sim b(u)$ for $\lim_{u\to\infty}
a(u)/b(u) =1$),
\begin{eqnarray}\label{FLP26}
&&[(\mathrm{e}^{-1}u)^{d-1}-(-1+\mathrm{e}^{-1}u)^{d-1}]\bigl(u\log(1+\log
|u|)\bigr)^{{1/2}}\nonumber
\\[-8pt]\\[-8pt]
&&\quad\sim(d-1)(\mathrm{e}^{-1})^{d-2}u^{d-{3/2}}(\log
\log(u))^{{1/2}} ,\nonumber
\end{eqnarray}
which ensures the existence of the two integrals in $(\ref{FLP24w})$.

Letting $t\to-\infty,$ we obtain \eqref{I1}. Next we calculate
\begin{eqnarray*}
&&\frac{1}{|t|^\alpha}\int_{t}^0(-s)^{d-1}|L(s)
|\,\mathrm{d}s=\frac{1}{|t|^\alpha}\int_{t}^T(-s)^{d-1}|L(s)
|\,\mathrm{d}s+\frac{1}{|t|^\alpha}\int_{T}^0(-s)^{d-1}|L(s)|\,\mathrm{d}s.
\end{eqnarray*}
The second term tends to zero as $t\to-\infty$, and we consider the first:
\begin{eqnarray*}
\frac{1}{|t|^\alpha}\int_{t}^T(-s)^{d-1}|L(s)|\,\mathrm{d}s&\leq&
\frac{M}{|t|^\alpha}\int_{t}^T(-s)^{d-1}(2|s|\log\log|s|)^{
{1/2}}\,\mathrm{d}s
\\
&\leq&\frac{M(2|t|\log\log|t|)^{{1/2}}}{|t|^\alpha}\int
_{t}^T(-s)^{d-1}\,\mathrm{d}s
\\
&=&\frac{M(2\log\log|t|)^{{1/2}}}{d|t|^{\alpha-(d+
{1/2})}}-\frac{|T|^dM(2\log\log|t|)^{{1/2}}}{d|t|^{\alpha
-{1/2}}}.
\end{eqnarray*}
Letting $t\to-\infty,$ we get \eqref{I2} and therefore the assertion.
\end{pf}

Theorem~\ref{th3.1} ensures the existence of the improper
Riemann--Stieltjes integral.

\begin{lemma}\label{le3.2}
Let $L^d$ be an FLP, $d\in(0,\frac{1}{2})$ and $\lambda>0$. Then for
$-\infty\leq a<\infty$
%
\begin{equation}\label{A4}
\int_a^t\mathrm{e}^{\lambda s}\,\mathrm{d}L^d_s ,\qquad t>a,
\end{equation}
exists a.s. as a Riemann--Stieltjes integral and is equal to
%
\begin{equation}\label{A5}
\mathrm{e}^{\lambda t}L^d_t-\mathrm{e}^{\lambda a}L^d_a-\lambda\int_a^tL^d_s\mathrm{e}^{\lambda
s}\,\mathrm{d}s .
\end{equation}
Furthermore, the function $(a,\infty)\to\mathbb{R}$ defined by
$t\mapsto\int_a^t\mathrm{e}^{\lambda s}\,\mathrm{d}L^d_s$ is continuous.
\end{lemma}

\begin{pf}
From Theorem~\ref{th3.1} we know that for all $\alpha>d+\frac{1}{2}$
there is a null set $N\subset\Omega$ such that for $\omega\in\Omega
\setminus N$ we have
%
\begin{equation}
\lim_{t\rightarrow-\infty}\frac{L^d_t(\omega)}{|t|^\alpha}=0
\end{equation}
and, hence, for all $\omega\in\Omega\setminus N$ and $t>a$, the
integral $\int_a^tL^d_u(\omega)\mathrm{e}^{\lambda u}\,\mathrm{d}u$ exists as a
Riemann--Stieltjes integral. For a compact interval $[a,t]$ this is
clear. Now consider $a=-\infty$. It suffices to show that $\int
_{-\infty}^TL^d_u(\omega)\mathrm{e}^{\lambda u}\,\mathrm{d}u$ exists for $T<-1$. This
follows from the inequality
\begin{eqnarray*}
\biggl|\int_{R}^TL^d_u(\omega)\mathrm{e}^{\lambda u}\,\mathrm{d}u\biggr|\leq
\int_{R}^T\frac{|L^d_u(\omega)|}{|u|^\alpha}\mathrm{e}^{\lambda u}|u|^\alpha \,\mathrm{d}u
\leq C\int_{R}^T\mathrm{e}^{\lambda u}|u|^\alpha \,\mathrm{d}u
\end{eqnarray*}
for some constant $C>0$, and the integral on the right-hand side exists
for $R\to-\infty$. Similarly,
%
\begin{equation}
\lim_{a\rightarrow-\infty}\mathrm{e}^{\lambda a}L^d_a(\omega)=0 .
\end{equation}
Now it follows by Wheeden and Zygmund \cite{whee}, Theorem 2.21, that
$(\ref{A4})$ also exists as a Riemann--Stieltjes integral and is equal
to (\ref{A5}). Since $(\ref{A5})$ is continuous in $t$ for all $t>a$,
the result is proven.
\end{pf}

Now we are ready to define the central object of this paper.
Recall that all integrals are Riemann--Stieltjes integrals based on
Lemma~\ref{le3.2}.

\begin{definition}\label{def3.3}
Let $L^d$ be an FLP, $d\in(0,\frac{1}{2})$ and $\lambda>0$. Then
%
\begin{equation}
\mathcal{L}^{d,\lambda}_t:=\int_{-\infty}^t\mathrm{e}^{-\lambda
(t-s)}\,\mathrm{d}L^d_s ,\qquad t\in\mathbb{R},
\end{equation}
is called \emph{an FLOUP}.
\end{definition}

Before returning to the Langevin equation in connection with the FLOUP,
we present some distributional properties of $\mathcal{L}^{d,\lambda
}$. With a little effort one can prove that $\mathcal{L}^{d,\lambda}$
is stationary, i.e., for all $t_1<\cdots<t_m$, $m\in\mathbb{N}$,
$h\in\mathbb{R}$,
%
\begin{equation}\label{le3.5}
(\mathcal{L}^{d,\lambda}_{t_1},\ldots,\mathcal{L}^{d,\lambda}_{t_n})
\stackrel{d}{=}(\mathcal{L}^{d,\lambda}_{t_1+h},\ldots,\mathcal
{L}^{d,\lambda}_{t_n+h}).
\end{equation}
For details see \cite{fink}, Lemma~6.1.3.

Although we mainly concentrate on Riemann--Stieltjes integrals, there
exist several results based on integrals in the $L^2(\Omega)$-sense
that we can use. Fractional integration can be considered as a
transformation of classical Riemann--Liouville fractional integrals,
which are defined for $0<\alpha<1$ by
\begin{eqnarray}
(I^\alpha_{-}f)(x)&=&\frac{1}{\Gamma(\alpha)}\int_x^\infty
f(t)(t-x)^{\alpha-1}\,\mathrm{d}t\nonumber
\quad \mbox{and}
\\[-8pt]\\[-8pt]
(I^\alpha_{+}f)(x)&=&\frac{1}{\Gamma(\alpha)}\int_{-\infty}^x
f(t)(t-x)^{\alpha-1}\,\mathrm{d}t,\nonumber
\end{eqnarray}
if the integrals exist for almost all $x\in\mathbb{R}$.
This is, for instance, the case if $f\in L^p(\mathbb{R})$ with $1\leq
p\leq\frac{1}{\alpha}$. The following result is a Riemann--Stieltjes
version of Theorem~3.5 of \cite{tina}.

\begin{proposition}\label{cor3.6}
Let $\mathcal{L}^{d,\lambda}$ be an FLOUP driven by an FLP $L^d$ of
bounded $p$-variation, $d\in(0,\frac{1}{2})$, $\lambda>0$ and $p>0$.
Then its finite-dimensional distributions have a characteristic function
\begin{eqnarray*}\label{FLP4}
&&E\Biggl[\exp\Biggl\{\sum_{j=1}^m\mathrm{i}u_j\mathcal{L}^{d,\lambda}_{t_j}
\Biggr\}\Biggr] =\exp\Biggl\{\int_\mathbb{R}
\psi_L\Biggl(\sum_{j=1}^mu_j\bigl(I^d_{-}\mathrm{e}^{-\lambda(t_j-\cdot)}1_{\{
t_j\geq\cdot\}}\bigr)(s)
\Biggr)\,\mathrm{d}s\Biggr\},
\\
&&\quad u_1,\ldots,u_m\in\mathbb{R},
\end{eqnarray*}
for $-\infty<t_1<\cdots<t_m<\infty$ and $\psi_L(u):=\int_{\mathbb
{R}}(\mathrm{e}^{\mathrm{i}ux}-1-\mathrm{i}ux)\nu(\mathrm{d}x)$, where $\nu$ is the L\'evy
measure of $L$.
Furthermore, for every $t\in\mathbb{R}$, the random variable
$\mathcal{L}^{d,\lambda}_t$ is infinitely divisible with a
characteristic triple given by $(\gamma_{\mathcal{L}^{d,\lambda
}}^t,0,\nu_{\mathcal{L}^{d,\lambda}}^t)$, where
%
\begin{eqnarray}\label{ctrip}
\gamma_{\mathcal{L}^{d,\lambda}}^t&=&-\int_\mathbb{R}\int_\mathbb
{R}\bigl(I^d_{-}\mathrm{e}^{-\lambda(t-\cdot)}1_{\{t\geq\cdot\}}\bigr)(s)x1_{\{
|(I^d_{-}\mathrm{e}^{-\lambda(t-\cdot)}1_{\{t\geq\cdot\}})(s)x|>1\}}\,\mathrm{d}\nu
(x)\,\mathrm{d}s ,\qquad 
\\
\nu_{\mathcal{L}^{d,\lambda}}^t(B)&=&\int_\mathbb{R}\int_\mathbb
{R}1_B\bigl(\bigl(I^d_{-}\mathrm{e}^{-\lambda(t-\cdot)}1_{\{t\geq\cdot\}}\bigr)(s)x\bigr)\,\mathrm{d}\nu
(x)\,\mathrm{d}s \quad
\forall\mbox{Borel sets $B$ in $\mathbb{R}$}.
\end{eqnarray}
\end{proposition}

\begin{pf}
For simple functions, the Riemann--Stieltjes integral and the
$L^2(\Omega)$-integral agree a.s. (see \cite{tina}, Proposition~5.2).
Now approximate the function $\mathrm{e}^{-\lambda(t-s)}1_{\{t\geq s\}}$ by
simple functions. While a.s. and $L^2(\Omega)$-convergence of the
Riemann--Stieltjes sums imply both convergence in probability, the
integrals are equal in probability and thus in distribution. Therefore
the result follows as in Theorem~3.5 of \cite{tina}.
\end{pf}

We now turn to the second-order properties of an FLOUP.
Cheridito, Kawaguchi and Maejima~\cite{cheri} present details
concerning the long memory property of an OU process driven by FBM.
Similarly, we shall show that the increments of the FLOUP exhibit
long-range dependence.
First, however, we need the following result (see also Proposition 4.4
of \cite{mat} and Proposition~5.7 of \cite{tina}).

\begin{theorem}\label{th3.11}
Let $L^d$ be an FLP of bounded $p$-variation, $d\in(0,\frac{1}{2})$,
$\lambda>0$, $p>0$ and $f,g:\mathbb{R}\to\mathbb{R}$ with
$|f|,|g|\in\mathfrak{W}^{\mathrm{con}}_q(\mathbb{R})$ for
$p^{-1}+q^{-1}>1$,
such that $\int_\mathbb{R}f(s)\,\mathrm{d}L^d_s$ and $\int_\mathbb
{R}g(s)\,\mathrm{d}L^d_s$ exist as Riemann--Stieltjes integrals.
Then we have
\begin{eqnarray}
&&E\biggl[\int_\mathbb{R}f(t)\,\mathrm{d}L^d_t\int_\mathbb{R}g(s)\,\mathrm{d}L^d_s
\biggr]\nonumber
\\[-8pt]\\[-8pt]
&&\quad=\frac{\Gamma(1-2d)E[(L(1))^2]}{\Gamma(d)\Gamma(1-d)}\int_\mathbb
{R}\int_\mathbb{R}f(t)g(s)|t-s|^{2d-1}\,\mathrm{d}s\,\mathrm{d}t .\nonumber
\end{eqnarray}
\end{theorem}

\begin{pf}
The proof follows again by using approximating simple functions and the
fact that
%
\begin{equation}
\int_{-\infty}^{t\wedge s}(t-u)^{d-1}(s-u)^{d-1}\,\mathrm{d}u =\frac{\Gamma
(d)\Gamma(1-2d)}{\Gamma(1-d)}|t-s|^{2d-1}
\end{equation}
for $t,s\in\mathbb{R}$, which can be found in Gripenberg and
Norros~\cite{gp}, page~404.
\end{pf}

Now we have everything together to derive the covariance structure of
an FLOUP.
The lengthy calculation works in a manner similar to that of Theorem
2.3 of \cite{cheri}. The asymptotic
is up to a multiplicative factor the same as in the FBM case.

\begin{theorem}\label{th3.12}
Let $L^d$ be an FLP, $d\in(0,\frac{1}{2})$, $\lambda>0$ and
$\mathcal{L}^{d,\lambda}$ the corresponding FLOUP. Then for $N\in
\mathbb{N}_0$ and for fixed $t\in\mathbb{R}$ we have as $ s\to
\infty$
%
\begin{eqnarray*}
&&\operatorname{Cov}(\mathcal{L}^{d,\lambda}_t,\mathcal
{L}^{d,\lambda}_{t+s})\nonumber
\\
&&\quad=
\frac{\Gamma(1-2d)E[L(1)]^2}{2d(2d+1)\Gamma(d)\Gamma(1-d)}
\sum_{n=1}^{N}\Biggl(\prod_{k=0}^{2n-1}(2d+1-k)\Biggr)\lambda^{-2n}
s^{2d+1-2n}
+\mathrm{O}(s^{2d-2N-1}).\nonumber
\end{eqnarray*}
\end{theorem}

Now it is clear that the increments of an FLOUP exhibit long-range
dependence in the sense of a non-summability property of the
autocovariance function.

We now return to the Langevin equation presented in \eqref{lange}.

\begin{theorem}\label{th3.8}
Let $L^d$ be an FLP, $d\in(0,\frac{1}{2})$ and $\lambda>0$. Then the
unique stationary pathwise solution of \eqref{lange}
is given a.s. by the corresponding FLOUP
\begin{eqnarray*}\label{A12}
\mathcal{L}^{d,\lambda}_t=\int_{-\infty}^t\mathrm{e}^{-\lambda(t-u)}\,\mathrm{d}L^d_u
,\qquad t\in\mathbb{R}.
\end{eqnarray*}
\end{theorem}

\begin{pf}
From Lemma~\ref{le3.2} we know that $\int_{-\infty}^t\mathrm{e}^{-\lambda
(t-u)}\,\mathrm{d}L^d_u$ exists for all $t\in\mathbb{R}$ a.s. as a
Riemann--Stieltjes integral. We fix $s\in\mathbb{R}$ and consider the
pathwise SDE
%
\begin{equation}\label{A13}
\mathcal{L}^{d,\lambda}_t=\xi_s-\lambda\int_s^t\mathcal
{L}^{d,\lambda}_u\,\mathrm{d}u+L^d_t-L^d_s,\qquad s\leq t,
\end{equation}
where $\xi_s:=\int_{-\infty}^s\mathrm{e}^{-\lambda(s-u)}\,\mathrm{d}L^d_u$. Obviously,
$\xi_s\in L^0(\Omega)$.
By arguments similar to those in the proof of \cite{cheri}, Proposition~A.1, we obtain
\begin{eqnarray*}
\mathcal{L}^{d,\lambda}_t=\mathrm{e}^{-\lambda t}\biggl\{\mathrm{e}^{\lambda s}\int
_{-\infty}^s\mathrm{e}^{-\lambda(s-u)}\,\mathrm{d}L^d_u+\int_{s}^t\mathrm{e}^{\lambda
u}\,\mathrm{d}L^d_u\biggr\}=\int_{-\infty}^t\mathrm{e}^{-\lambda(t-u)}\,\mathrm{d}L^d_u ,\qquad
t\in\mathbb{R},
\end{eqnarray*}
is the unique pathwise solution of (\ref{A13}) and, therefore, by
\eqref{le3.5} a stationary solution of (\ref{lange}).

On the other hand, let $(X_t)_{t\in\mathbb{R}}$ be a stationary
solution of (\ref{lange}). We show that $(X_t)_{t\in\mathbb
{R}}=(\mathcal{L}^{d,\lambda}_t)_{t\in\mathbb{R}}$ holds for almost
all $\omega\in\Omega$.
Set $A:=\{\omega\in\Omega\dvtx (X_t(\omega))_{t\in\mathbb{R}}\neq
(\mathcal{L}^{d,\lambda}_t(\omega))_{t\in\mathbb{R}}\}$ and assume
that $P(A)>0$. For $\omega\in A$ fix $t\in\mathbb{R}$ with
$X_t(\omega)\neq\mathcal{L}^{d,\lambda}_t(\omega)$.
Then we have for $s\leq t$ by (\ref{A13})
\begin{eqnarray*}
0&\neq&|X_t-\mathcal{L}^{d,\lambda}_t|=
\biggl|\mathrm{e}^{-\lambda t}\biggl\{\mathrm{e}^{\lambda s} X_s+\int_s^t\mathrm{e}^{\lambda
u}\,\mathrm{d}L^d_u\biggr\}-\int_{-\infty}^t\mathrm{e}^{-\lambda(t-v)}\,\mathrm{d}L^d_v\biggr|
\\
&=&\biggl|\mathrm{e}^{-\lambda(t-s)}X_s-\int_{-\infty}^s\mathrm{e}^{-\lambda
(t-u)}\,\mathrm{d}L^d_u\biggr|
=\mathrm{e}^{-\lambda t}\mathrm{e}^{\lambda s}|X_s-\mathcal{L}^{d,\lambda
}_s|,
\end{eqnarray*}
where we supressed the chosen $\omega$ for simplicity.
Since $\lambda>0$ and $s\to-\infty$ we conclude that $|X_s(\omega
)-\mathcal{L}^{d,\lambda}_s(\omega)|\to\infty$ for $s\to-\infty
$. Therefore, on $A$ we have $|X_t-\mathcal{L}^{d,\lambda}_t|\to
\infty$ for $t\to-\infty$. 
For a given $K>0$ we define $\omega$-wise the random variable $T\dvtx A\to
\mathbb{R}$ with $|X_t-\mathcal{L}^{d,\lambda}_t|\geq\frac
{K}{P(A)}$ for $t\leq T$ on $A$. Hence,
\begin{eqnarray*}\label{qq1}
\nonumber E|X_t-\mathcal{L}^{d,\lambda}_t|&\geq&E
\bigl[|X_t-\mathcal{L}^{d,\lambda}_t|1_{\{t\leq T\}}1_A\bigr]+E
\bigl[|X_t-\mathcal{L}^{d,\lambda}_t|1_{\{t>T\}}1_A\bigr]\geq\frac
{K}{P(A)}P(\{t\leq T\}\cap A).
\end{eqnarray*}
Furthermore, we know that $\{t\leq T\}\cap A\subseteq\{s\leq T\}\cap
A$ for $s\leq t$. Choosing a sequence $(t_n)_{n\in\mathbb{N}}$ of
real numbers with $\lim_{n\to\infty}t_n=-\infty$ we get by
continuity of $P$
\begin{eqnarray*}
\lim_{n\to\infty} P(\{t_n\leq T\}\cap A) =P\biggl(\bigcup_{n\in
\mathbb{N}}\{t_n\leq T\}\cap A\biggr)=P(A).
\end{eqnarray*}
Putting everything together we arrive at
%
\begin{eqnarray}
\lim_{n\to\infty} E|X_{t_n}-\mathcal{L}^{d,\lambda}_{t_n}|\geq
\lim_{n\to\infty} \frac{K}{P(A)}P(\{t_n\leq T\}\cap A)=K.
\end{eqnarray}
Hence, $\lim_{n\to\infty} E|X_{t_n}-\mathcal{L}^{d,\lambda
}_{t_n}|=\infty$. However, we have now
\begin{eqnarray*}
E|X_{t_n}|=E|X_{t_n}-\mathcal{L}^{d,\lambda}_{t_n}-(-\mathcal
{L}^{d,\lambda}_{t_n})|\geq E|X_{t_n}-\mathcal{L}^{d,\lambda
}_{t_n}|-E|\mathcal{L}^{d,\lambda}_{t_n}|,
\end{eqnarray*}
where $E|\mathcal{L}^{d,\lambda}_{t_n}|$ is independent of $t_n$.
Thus, $\lim_{n\to\infty}E|X_{t_n}|=\infty$ and, by stationarity,
$E|X_t|=\infty$ for all $t\in\mathbb{R}$.
However, we also have for fixed $s\leq t$
%
\begin{eqnarray}
\lim_{t\to\infty}(X_t-\mathcal{L}^{d,\lambda}_t)=\lim_{t\to
\infty}\mathrm{e}^{-\lambda t}\biggl(\mathrm{e}^{\lambda s}X_s-\int_{-\infty
}^s\mathrm{e}^{\lambda u}\,\mathrm{d}L^d_u\biggr)=0\qquad\mbox{a.s.}
\end{eqnarray}
Hence, by stationary $X_t\stackrel{d}{=}\mathcal{L}^{d,\lambda}_t$,
but $E|\mathcal{L}^{d,\lambda}_t|<\infty$, which is a contradiction
and, thus, we conclude that $P(A)=0$.
\end{pf}

The following Ornstein--Uhlenbeck operator will be used in the next
section to obtain solutions to SDEs with different starting values. The
operator here modifies the starting value of the FLOUP and the lemma
shows that this modified process still solves the Langevin equation.

\begin{definition}\label{def3.9}
Let $L^d$ be an FLP, $d\in(0,\frac{1}{2})$, $\lambda>0$ and
$\mathcal{L}^{d,\lambda}$ the corresponding FLOUP.
We define the Ornstein--Uhlenbeck operator by
\begin{eqnarray}\label{OUO1}
&&\mathfrak{L}^\lambda(L^d,\cdot,\cdot)\dvtx \mathbb{R}\times\mathbb
{R}\longrightarrow \mathcal{C}^0(\mathbb{R}),\nonumber
\\[-8pt]\\[-8pt]
&&(\tau,z)\longmapsto\mathcal{L}^{d,\lambda}_t-\mathrm{e}^{-\lambda(t-\tau)}
\mathcal{L}^{d,\lambda}_\tau+\mathrm{e}^{-\lambda(t-\tau)}z.\nonumber
\end{eqnarray}
\end{definition}

It is immediate from this definition that $\mathfrak{L}_\tau^\lambda
(L^d,\tau,z)=z$ a.s. for $(\tau,z)\in\mathbb{R}^2$.

The next lemma shows that $L^d$ transformed by the Ornstein--Uhlenbeck
operator still satisfies the Langevin equation; its proof follows
directly by the definition.

\begin{lemma}\label{le3.10}
Let $L^d$ be an FLP, $d\in(0,\frac{1}{2})$, $\lambda>0$ and
$\mathcal{L}^{d,\lambda}$ be the corresponding FLOUP. For a
continuous process $l:=(l_t)_{t\in\mathbb{R}}$ the identity
$l_t=\mathfrak{L}^\lambda_t(L^d,\tau,l_\tau)$ holds for all $\tau
,t\in\mathbb{R}$ if and only if
%
\begin{equation}\label{SDE2}
\mathrm{d}l_t=-\lambda l_t\,\mathrm{d}t+\mathrm{d}L^d_t ,\qquad t\in\mathbb{R}.
\end{equation}
\end{lemma}

\section{Solutions of fractional SDEs by state space transforms and
proper triples}\label{sec4}

In this section we start with SDEs driven by FLPs.
Using pathwise integration we must solve for almost all $\omega\in
\Omega$ a deterministic integral equation.
Consequently, we build on an extensive theory starting with the seminal
work by Young~\cite{Young}. We also recall that for Brownian motion
the pathwise approach goes back to \cite{doss} and \cite
{sussman} leading to the Fisk--Stratonovich integral.
Required is that $\mu$ is Lipschitz-continuous and ${\sigma}\in
C^2(\mathbb{R}
)$ with bounded first and second derivatives.
Readable accounts on the history can be found in \cite{KS} and in Ikeda and \cite{IKWA}.

Regularity assumptions of sample paths of the driving process like H\"
older continuity or bounded $p$-variation for $p<2$ have been
considered by
Lyons~\cite{Lyons}.
We shall work in the framework of $p$-variation, however, to go beyond
the work of Lyons, who proves only existence and uniqueness theorems
under certain Lipschitz assumptions on the coefficient functions and
gives no analytical form of the solution.

The approach by Z\"ahle \cite{zaehle2} is indeed comparable to ours,
where explicit solutions can be given under differentiability and
Lipschitz conditions on the coefficient functions.
Most of her results can be applied to SDEs driven by an FLP of bounded
$p$-variation for $p<2$.
We believe that the contribution of our work is two-fold. First, our
assumptions are easy to verify and, second, we are able to present
analytic solutions to SDEs of the form
%
\begin{equation}\label{FLP1211}
\mathrm{d}X_t=(\alpha|X_t|^\gamma+\beta X_t)\,\mathrm{d}t+\sigma|X_t|^\gamma\, \mathrm{d}L^d_t
,\qquad t\in\mathbb{R}.
\end{equation}
In this situation we cannot apply the results of \cite{zaehle2},
since the volatility coefficient does not match the required
differentiability assumption.
Lyons \cite{Lyons} provides us at least with an existence theorem, but
gives no closed form solution.

Aiming at solutions to similar SDEs, driven however by FBM, Buchmann
and Kl\"uppelberg~\cite{KL} presented a theory that can be modified to
cover SDEs driven by FLPs.
The idea is to present solutions to SDEs like, for instance, \eqref
{FLP1211} as monotone transformations of the FLOUP.
The question we shall answer is, given an SDE
%
\begin{eqnarray}\label{sdeflouptrans}
\mathrm{d}X_t = \mu(X_t) + {\sigma}(X_t) \,\mathrm{d}L_t^d
\end{eqnarray}
for specific $\mu$ and $\sigma$,
which monotone transformation of the FLOUP is a solution to \eqref
{sdeflouptrans}?

First we have to establish certain regularity conditions on the
coefficients $\mu$ and $\sigma$.

\begin{definition}\label{def4.1}
\textup{(i)}   A triple $(I,\mu,\sigma)$ is called \emph{strongly proper} if
and only if it satisfies the following properties:
\begin{enumerate}[(P1)]
\item[(P1)] $I=(a,b)\subseteq\mathbb{R}$ is an open interval, where
$-\infty\leq a<b\leq\infty$ and $\mu,\sigma\in\mathcal{C}^0(I)$.
\item[(P2)] There exists a strictly decreasing $\psi$, absolutely
continuous with respect to the Lebesgue measure such that $\psi=\mu
/\sigma$ on $I\setminus Z(\sigma)$ where $Z(\sigma)$ are the zeros
of $\sigma$, and
\begin{eqnarray*}
\lim_{x\nearrow b}\psi(x)=-\lim_{x\searrow a}\psi(x)=-\infty .
\end{eqnarray*}
\item[(P3)]There exists $\lambda>0$ such that $\sigma\psi'\equiv
\lambda$ holds on $I$ Lebesgue-a.e.
\item[(P4)]The inverse function $\psi^{-1}\dvtx \mathbb{R}\rightarrow
\psi^{-1}(\mathbb{R})=I$ is differentiable and $(\psi^{-1})'\in
\operatorname{Lip}(\mathbb{R})$.
\end{enumerate}

\hspace*{2,3pt}\textup{(ii)}   We call the triple $(I,\mu,\sigma)$ \emph{proper} if only
\textup{(P1)--(P3)} are satisfied.

\textup{(iii)}   The interval $I$ is called the \emph{state space}, the unique
constant $\lambda>0$ in \textup{(P3)} is called the {friction coefficient
(FC)} and the unique function $f\dvtx \mathbb{R}\rightarrow I=f(\mathbb
{R})$, $f(x):=\psi^{-1}(-\lambda x)$, is called the {state space
transform (SST)} for $(I,\mu,\sigma)$.
\end{definition}

Condition (P4) differs from the $H$-proper assumption required in
\cite{KL}, because we work with
$p$-variation instead of H\"older continuity.

As pointed out in \cite{KL},
$\psi\dvtx I\to\psi(I)=\mathbb{R}$ is by (P2) strictly decreasing and
a.e. differentiable on $I$ with $\psi'\leq0$.
Condition (P3) implies that $Z(\sigma)$ and $Z(\psi')$ have Lebesgue
measure zero.
Also we have that $\sigma$ is non-negative and $1/\sigma\in\mathcal
{L}_C(I)$, where $\mathcal{L}_C(I)$ denotes the locally integrable
functions on $I$; $I\setminus Z(\sigma)$ is dense and open in $I$ by (P1).
Therefore, the equality $\mu=\sigma\psi$ extends to $I$.
It follows that $\psi$ and $\lambda$ are uniquely determined by $\mu
$ and $\sigma$.

As can be seen from Definition~\ref{def4.1} our coefficient functions
are only defined on the interval~$I$, which can be any interval in
$\mathbb{R}$. To account for this situation we need to specify what
will be understood to be a solution to an SDE.

\begin{definition}\label{def4.2}
Let $L^d$ be an FLP of bounded $p$-variation, $p\in[1,2)$ and $d\in
(0,\frac{1}{2})$.
Suppose that $I\subseteq\mathbb{R}$ is a non-empty interval and $\mu
,\sigma\in\mathcal{C}^0(I)$. We refer to a stochastic process
$X:=(X_t)_{t\in\mathbb{R}}$ as a pathwise solution of the SDE
%
\begin{equation}\label{FP111}
\mathrm{d}X_t=\mu(X_t)\,\mathrm{d}t+\sigma(X_t)\,\mathrm{d}L^d_t,\qquad t\in\mathbb{R},
\end{equation}
if for almost all sample paths the following conditions are satisfied:
$X\in\mathfrak{W}^{\mathrm{con}}_p(\mathbb{R})$ and the image of
$X$ is a subset of $I$ such that for $s\leq t$:
\begin{enumerate}[(S1)]
\item[(S1)] $\sigma\circ X$ is a.s. Riemann--Stieltjes integrable
with respect to $L^d$ on $[s,t]$;
\item[(S2)] The following integral equation holds:
\begin{eqnarray*}
X_t-X_s=\int_s^t\mu(X_u)\,\mathrm{d}u+\int_s^t\sigma(X_u)\,\mathrm{d}L^d_u.
\end{eqnarray*}
\end{enumerate}
The space of all solutions of \textup{(\ref{FP111})} is denoted by $\mathcal
{S}(I,\mu,\sigma,L^d)$.
\end{definition}

We consider now an SDE as given in $(\ref{FP111})$. If we assume that
the triple $(I,\mu,\sigma)$ is strongly proper with SST $f$ and FC
$\lambda$, we define the following operator
%
\begin{eqnarray}\label{FF1}
X^{f,\lambda}(L^d,\cdot,\cdot)\dvtx \mathbb{R}\times I \longrightarrow
\mathcal{C}^0(\mathbb{R}),\qquad(\tau,z)\longmapsto f(\mathfrak
{L}_t^\lambda(
L^d,\tau,f^{-1}(z)))
\end{eqnarray}
with Ornstein--Uhlenbeck operator $\mathfrak{L}_t^\lambda$ as in
Definition~\ref{def3.3}.
We also remark that
%
\begin{equation}\label{FF2}
X_t^{f,\lambda}(L^d,\tau,f(\mathcal{L}_\tau^{d,\lambda
}))=f(\mathcal{L}^{d,\lambda}_t),\qquad t\in\mathbb{R}.
\end{equation}
Before we state our main results we state the following technical
lemma, which will be needed in the proofs.

\begin{lemma}[(Version of Lemma 3.2 \cite{KL})]\label{le4.5}
Let $(I,\mu,\sigma)$ be a strongly proper triple with the
corresponding SST $f$. Then
$f\in\mathcal{C}^1(\mathbb{R})$ with derivative $f'=\sigma\circ f$.
Also $f^{-1}\in\mathcal{C}^1(I\setminus Z(\sigma))$ with
$(f^{-1})'(x)=1/\sigma(x)$ for all $x\in I\setminus Z(\sigma)$.
\end{lemma}

Next we state the existence theorem. Let $\mathrm{M}(\Omega,I)$
denote all mappings from $\Omega$ into $I$.

\begin{theorem}\label{th4.6}
Let $L^d$ be an FLP of bounded $p$-variation, $p\in[1,2)$ and $d\in
(0,\frac{1}{2})$. If $(I,\mu,\sigma)$ is strongly proper with SST
$f$ and FC $\lambda>0$, then
\begin{eqnarray*}\label{FF3}
\{X^{f,\lambda}(L^d,\tau,W)\dvtx \tau\in\mathbb{R},W\in
{M}(\Omega,I)\}\subseteq\mathcal{S}(I,\mu,\sigma,L^d).
\end{eqnarray*}
\end{theorem}

\begin{pf} Because we consider pathwise solutions we can w.l.o.g. assume
that $W=z$ a.s. for some $z\in I$.
Now fix $\tau\in\mathbb{R}$ and $z\in I$ and define
\begin{eqnarray*}
l_t:=\mathfrak{L}_t^\lambda(L^d,\tau,f^{-1}(z))\quad\mbox{and}\quad
Y_t:=X_t^{f,\lambda}(L^d,\tau,z),\qquad t\in\mathbb{R}.
\end{eqnarray*}
We show that $Y=(Y_t)_{t\in\mathbb{R}}\in\mathcal{S}(I,\mu,\sigma
,L^d)$. Obviously $Y$ takes only values in $I$.
Since $f\in\mathcal{C}^1(\mathbb{R})$ and $l$ is of bounded
$p$-variation, we know that $Y=f\circ l\in\mathfrak{W}_p^{\mathrm
{con}}(\mathbb{R})$.
With the chain rule from Theorem~\ref{thA.3} we get for $s,t\in
\mathbb{R}$
%
\begin{equation}\label{FF6}
Y_t-Y_s=f(l_t)-f(l_s)=\int_s^t f'(l_u)\,\mathrm{d}l_u ,\qquad s\leq t,
\end{equation}
since $l$ solves $(\ref{SDE2})$,
%
\begin{equation}\label{FF7}
l_u=l_s-\lambda\int_s^ul_v\,\mathrm{d}v+L^d_u-L^d_s,\qquad s\leq u.
\end{equation}
The Riemann--Stieltjes integral is additive with respect to a sum of
integrators, if the Riemann--Stieltjes integrals exist separately for
each integrator.
This is true in our case, because $\int_s^ul_v\,\mathrm{d}v$ is of finite
variation and $L^d_u$ is of bounded $p$-variation. Since $f'(l_u)$ is
continuous and also of bounded $p$-variation, $(\ref{FF6})$ and $(\ref
{FF7})$ imply for $s,t\in\mathbb{R}$
\begin{eqnarray*}
Y_t-Y_s=\int_s^t f'(l_u)\,\mathrm{d}\biggl(-\lambda\int_s^ul_v\,\mathrm{d}v\biggr)+\int
_s^t f'(l_u)\,\mathrm{d}L^d_u ,\qquad s\leq t.
\end{eqnarray*}
Furthermore, $\int_s^ul_v\,\mathrm{d}v$ is differentiable and $f'(l_u)l_u$ is
continuous as a function of $u$ and, thus, we get by the density
formula for Riemann--Stieltjes integrals on compacts for $s,t\in
\mathbb{R}$
\begin{eqnarray*}
Y_t-Y_s=-\lambda\int_s^t f'(l_u)l_u\,\mathrm{d}u+\int_s^t f'(l_u)\,\mathrm{d}L^d_u ,\qquad
s\leq t.
\end{eqnarray*}
From Lemma~\ref{le4.5} we obtain $f'=\sigma\circ f$, hence $\sigma
\circ Y=f'\circ l\in\mathfrak{W}_p^{\mathrm{con}}(\mathbb{R})$.
By Definition~\ref{def4.1}(P3) and the interpretation following this
definition we find that $\sigma f^{-1}=-\sigma\psi/\lambda=-\mu
/\lambda$.
This yields
\begin{eqnarray*}\label{FF10}
\nonumber Y_t-Y_s&=&-\lambda\int_s^t f'(l_u)l_u\,\mathrm{d}u+\int_s^t
f'(l_u)\,\mathrm{d}L^d_u
\\ \nonumber
&=&-\lambda\int_s^t\sigma(f(l_u))l_u\,\mathrm{d}u+\int_s^t \sigma
(f(l_u))\,\mathrm{d}L^d_u
\\
&=&\int_s^t \mu(Y_u)\,\mathrm{d}u+\int_s^t \sigma(Y_u)\,\mathrm{d}L^d_u ,\qquad s\leq t,
\end{eqnarray*}
where we used in the last line the equality $\sigma(f(l_u))l_u=\sigma
(Y_u)f^{-1}(Y_u)=-\mu(Y_u)/\lambda$.
Finally, we have $Y\in\mathcal{S}(I,\mu,\sigma,L^d)$.
\end{pf}

The following result ensures uniqueness under natural conditions.

\begin{theorem}\label{th4.7}
Let $L^d$ be an FLP of bounded $p$-variation, $p\in[1,2)$ and $d\in
(0,\frac{1}{2})$. Let also $(I,\mu,\sigma)$ be strongly proper with
SST $f$ and FC $\lambda>0$. Furthermore, assume that $Z(\sigma
)=\varnothing$.
Then
\begin{eqnarray*}\label{FF3i}
\{X^{f,\lambda}(L^d,\tau,W)\dvtx \tau\in\mathbb{R},W\in
{M}(\Omega,I)\}= \mathcal{S}(I,\mu,\sigma,L^d).
\end{eqnarray*}
\end{theorem}

\begin{pf}
From $Z(\sigma)=\varnothing$, we know by Lemma~\ref{le4.5} that $f\in
\mathcal{C}^1(\mathbb{R})$ and $(f^{-1})'(x)=1/\sigma(x)$ for all
$x\in I$. Let $X\in\mathcal{S}(I,\mu,\sigma,L^d)$.
From Definition~\ref{def4.1} we know that $X\in\mathfrak
{W}_p^{\mathrm{con}}(\mathbb{R})$ a.s. and from $(f^{-1})'\in\operatorname
{Lip}(I)$ we get by the chain rule from Theorem~\ref{thA.3} for
$s,t\in\mathbb{R}$
%
\begin{equation}\label{FFFP1}
f^{-1}(X_t)-f^{-1}(X_s)=\int_s^tf^{-1}(X_u)\,\mathrm{d}X_u=\int_s^t\frac
{1}{\sigma(X_u)}\,\mathrm{d}X_u ,\qquad s\leq t.
\end{equation}
Since $X\in\mathcal{S}(I,\mu,\sigma,L^d)$, we know that for $s,u\in
\mathbb{R}$
%
\begin{equation}\label{FFFP2}
X_u=X_s+\int_s^u\mu(X_v)\,\mathrm{d}v+\int_s^u\sigma(X_v)\,\mathrm{d}L^d_v ,\qquad s\leq u.
\end{equation}
Now $\int_s^u\mu(X_v)\,\mathrm{d}v$ is of finite variation and by the density
formula of Theorem~\ref{thA.4} the integral $\int_s^u\sigma
(X_v)\,\mathrm{d}L^d_v$ is of bounded $p$-variation as a function of $u$. By
putting $(\ref{FFFP1})$ and $(\ref{FFFP2})$ together and using again
Theorem~\ref{thA.4}, we get for $s,t\in\mathbb{R}$
\begin{eqnarray*}
&&f^{-1}(X_t)-f^{-1}(X_s)\nonumber
\\
&&\quad =\int_s^t\frac{1}{\sigma(X_u)}\,\mathrm{d}
\biggl(X_s+\int_s^u\mu(X_v)\,\mathrm{d}v+\int_s^u\sigma(X_v)\,\mathrm{d}L^d_v\biggr)
\\
&&\quad=\int_s^t\frac{1}{\sigma(X_u)}\,\mathrm{d}\biggl(\int_s^u\mu
(X_v)\,\mathrm{d}v\biggr)+\int_s^t\frac{1}{\sigma(X_u)}\,\mathrm{d}\biggl(
\int_s^u\sigma(X_v)\,\mathrm{d}L^d_v\biggr)
\\
&&\quad=\int_s^t\frac{\mu(X_u)}{\sigma(X_u)}\,\mathrm{d}u+L^d_t-L^d_s,\qquad s\leq t,
\end{eqnarray*}
since $(I,\mu,\sigma)$ is proper, $\psi(x)=\mu(x)(\sigma(x))^{-1}$
and $\psi(x)=-\lambda f^{-1}(x)$ hold for all $x\in I$. Thus,
\begin{eqnarray*}
f^{-1}(X_t)-f^{-1}(X_s)&=&\int_s^t\frac{\mu(X_u)}{\sigma
(X_u)}\,\mathrm{d}u+L^d_t-L^d_s\nonumber
\\
\nonumber&=&\int_s^t\psi(X_u)\,\mathrm{d}u+L^d_t-L^d_s
\\
&=&-\lambda\int_s^tf^{-1}(X_u)\,\mathrm{d}u+L^d_t-L^d_s,\qquad s\leq t.
\end{eqnarray*}
Hence, $f^{-1}(X)$ is a solution of (\ref{SDE2}). Fixing $\tau\in
\mathbb{R}$ we see by Lemma~\ref{le3.10} that $f^{-1}(X)=\mathfrak
{L}^\lambda(L^d,\tau,f^{-1}(X_\tau))$ and, finally, $X=
X^{f,\lambda}(L^d,\tau,X_\tau)$.
\end{pf}

The next corollary covers the important case of a stationary solution.

\begin{corollary}\label{cor4.8}
Let $L^d$ be an FLP of bounded $p$-variation, $p\in[1,2)$ and $d\in
(0,\frac{1}{2})$ and $\mathcal{L}^{d,\lambda}$ be the corresponding FLOUP.
Furthermore, let $(I,\mu,\sigma)$ be strongly proper with SST $f$ and
FC $\lambda>0$.
Set $X_t=f(\mathcal{L}^{d,\lambda}_t)$ for $t\in\mathbb{R}$.
Then the following assertions hold:
\begin{enumerate}[(ii)]
\item[(i)] $X$ is a stationary pathwise solution of the SDE
%
\begin{equation}\label{wesd}
\mathrm{d}X_t=\mu(X_t)\,\mathrm{d}t+\sigma(X_t)\,\mathrm{d}L^d_t ,\qquad t\in\mathbb{R}.
\end{equation}
\item[(ii)] If $Z(\sigma)=\varnothing$, then $X$ is the unique
stationary pathwise solution of $(\ref{wesd})$.
\end{enumerate}
\end{corollary}

\begin{pf}
(i) From Theorem~\ref{th4.6} we know that $X=X^{f,\lambda
}(L^d,0,f(\mathcal{L}^{d,\lambda}_0))=f(\mathcal{L}^{d,\lambda})$
is a pathwise solution of $(\ref{wesd})$. Furthermore, $X$ is
stationary as a transformation of a stationary FLOUP.

(ii) Given a pathwise solution of $(\ref{wesd})$, Theorem~\ref{th4.7}
supplies us with a $W\in{M}(\Omega,I)$ such that
$Y_t=f(\mathrm{e}^{-\lambda t}W+\mathcal{L}^{d,\lambda}_t)$. If we want $Y$ to
be strictly stationary, we must have $W=0$ a.s. and get $Y=X$.
\end{pf}

\section{Examples of SDEs driven by FLPs}\label{sec5}

\subsection{Examples by means of strongly proper triples}\label{ss51}

This section is dedicated to examples, which we illustrate by
simulations. For those we consider as driving L\'evy process a
compensated Poisson process $L^{\theta}
$ with intensity $\theta>0$; that is,
\begin{eqnarray*}L^{\theta}_t
:=P^{\theta}_t-t\theta,\qquad t\in\mathbb{R},
\end{eqnarray*}
where $P^{\theta}$ is a L\'evy process with drift $\gamma=0$ and L\'evy
 measure $\nu(\mathrm{d}x)=\theta\delta_1(\mathrm{d}x)$ without Brownian component.
Of course, we consider this process to be defined on the whole of
$\mathbb{R}$ using $(\ref{qwas1})$.

In a first step we simulate sample paths of $L^\theta$ and compute the
corresponding FLP $L^d$ by a Riemann--Stieltjes approximation; that is,
we approximate
\begin{eqnarray*}
L_t^d&\approx&\frac{1}{\Gamma(d+1)}\Biggl\{\sum_{k=-n^2}^0
\biggl[\biggl(t-\frac{k}{n}\biggr)^d-\biggl(-
\frac{k}{n}\biggr)^d\biggr]\bigl(L^{a,b}_{
{(k+1)/n}}-L^{a,b}_{{k/n}}\bigr)
\\
&&{}\hspace*{33pt}\quad+\sum_{k=1}^{[nt]}\biggl(t-\frac{k}{n}
\biggr)^d\bigl(L^{a,b}_{{(k+1)/n}}
-L^{a,b}_{{k/n}}\bigr)\Biggr\},\qquad t\in\mathbb{R}.
\end{eqnarray*}
From Theorem~2.55 of \cite{tina} we know that
the quality of this approximation is
\begin{eqnarray*}
\mathrm{O}(n^{d-{1/2}})+\mathrm{O}(n^{-{1/2}})+\mathrm{O}\bigl(n^{{(1+2d-2d^2)/(2d-3)}}\bigr).
\end{eqnarray*}
Furthermore, the Poisson-FLP is of finite variation by Theorem~2.25 of
\cite{tina}.

Now we use a version of the explicit Euler method for the SDE \eqref{lange}
\begin{eqnarray*}
\mathrm{d}\mathcal{L}^{d,\lambda}_t=-\lambda\mathcal{L}^{d,\lambda}_t \,\mathrm{d}t +
\,\mathrm{d}L^d_t,\qquad t\in\mathbb{R},
\end{eqnarray*}
to compute sample paths of the FLOUP.
We want to remark that all these computations are pathwise. Probability
comes in only through the underlying paths of the driving L\'evy processes.

Next we study some examples of solutions to the SDE \eqref{FP111}
given by strongly proper triples.
We will mainly draw from structural results of \cite{KL} taking into account that their $H$-proper condition
has to be replaced by our assumption (P4) in Definition~\ref{def4.1}.

For the rest of this section, let $L^d$ be an FLP of bounded
$p$-variation, $p\in[1,2)$ and $d\in(0,\frac{1}{2})$.

\begin{example}\label{exe1}
As a first example, we consider for parameters ${\alpha},\beta\in
\mathbb{R}$
and ${\sigma}>0$ an SDE of the form
\begin{eqnarray*}
\mathrm{d}X_t=(\alpha|X_t|^\gamma+\beta X_t)\,\mathrm{d}t+\sigma|X_t|^\gamma \,\mathrm{d}L^d_t
,\qquad t\in\mathbb{R} .
\end{eqnarray*}
We analyse this SDE by taking the volatility coefficient $\sigma
\dvtx\mathbb{R}\to[0,\infty)$ defined by $\sigma(x):=\sigma
_0|x|^\gamma$ for $\sigma_0>0$ and $\gamma\in\mathbb{R}$ as given.
The question is now, what drift functions $\mu$ and intervals $I$ lead
to strongly proper triples $(I,\mu,\sigma)$ as defined in
Definition~\ref{def4.1}? More precisely, we want to find elements in
the set
\begin{eqnarray*}
\mathcal{K}^I_\sigma:=\{(\lambda,\mu)\in\mathbb{R}^+\times
\mathcal{C}^0(I)\dvtx (I,\mu,\sigma)\mbox{ is proper with FC $\lambda
$}\}.
\end{eqnarray*}
Using Proposition~5.5 of~\cite{KL} we see that only $\gamma\in[0,1]$
leads to a non-empty $\mathcal{K}^I_\sigma$.
We consider the cases $\gamma=0$, $\gamma=1$ and $\gamma\in(0,1)$
separately.

Take first $\gamma=0$. For a triple $(I,\mu,\sigma)$ to be proper we
must have that $\sigma\psi'\equiv-\lambda$ with $\psi=\mu/\sigma
$. This results in
$\mathrm{d}X_t=(\alpha+\beta X_t) \,\mathrm{d}t +{\sigma}\,\mathrm{d}L^d_t$ with state space
$I=\mathbb{R}$
and SST is affine, more precisely, $f(x)= \alpha+\beta x$ for $x\in
\mathbb{R}$.

If $\gamma=1$, by Proposition~5.6 of~\cite{KL}, the state space can
only be either $I=(-\infty,0)$ or $I=(0,+\infty)$.
For $I=(0,\infty)$ we get
\begin{eqnarray*}
\mathcal{K}^{(0,\infty)}_\sigma&=&\{(|\beta|,\mu)\in\mathbb
{R}^+\times\mathcal{C}^0(I)\dvtx \mu(x)=\alpha x+\beta x\log x,\alpha\in
\mathbb{R},\beta<0,x\in(0,\infty)\}
\end{eqnarray*}
and the state space transform is $f(x)=\exp\{\sigma_0x-\frac
{\alpha}{\beta}\}$ for $x\in(0,\infty)$.
Simple calculation ensures that condition (P4) of Definition~\ref
{def4.1} is satisfied, and every element of $\mathcal{K}^{(0,\infty
)}_\sigma$ leads to a strongly proper triple.
An example of an SDE of this kind for ${\alpha}=0$ can be found in
\eqref{eq5.6}.
The case $I=(-\infty,0)$ can be treated analogously.

Finally, we consider $\gamma\in(0,1)$. Proposition~5.8 of~\cite{KL}
shows that the only possible state space is the whole real line
$\mathbb{R}$ and
\begin{eqnarray*}
\mathcal{K}^{\mathbb{R}}_\sigma=\bigl\{\bigl((1-\gamma)|\beta|,\mu
\bigr)\in\mathbb{R}^+\times\mathcal{C}^0(I)\dvtx \mu(x)=\alpha|x|^\gamma
+\beta x,\alpha\in\mathbb{R},\beta<0,x\in\mathbb{R}\bigr\}.
\end{eqnarray*}
 Furthermore, the SST is given by
\begin{eqnarray*}
f(x)=\mathrm{sign}\biggl((1-\gamma)\sigma_0 x-\frac{\alpha}{\beta
}\biggr)\biggl|(1-\gamma)\sigma_0 x-\frac{\alpha}{\beta}
\biggr|^{{1/(1-\gamma)}}.
\end{eqnarray*}
The derivative of $f$ can easily be calculated yielding that only
$\gamma\in[\frac{1}{2},1)$ leads to strongly proper triples.
An example for such an SDE (with ${\alpha}=0$) is a fractional
Cox--Ingersoll--Ross-type model, which is investigated in detail in
Section~\ref{ss52} (cf. the SDE in \eqref{eq5.3}).
\end{example}

\begin{example}\label{exe2}
We consider the following SDEs with affine drift.
\begin{eqnarray*}
\mathrm{d}X_t=(\alpha+\beta X_t)\,\mathrm{d}t+\sigma(X_t)\,\mathrm{d}L^d_t ,\qquad t\in\mathbb{R} ,
\end{eqnarray*}
that is, $\mu\dvtx \mathbb{R}\to\mathbb{R}$ is defined by $\mu
(x):=\alpha+\beta x$ for $\alpha,\beta\in\mathbb{R}$.
To find suitable volatility coefficients and state spaces, we consider
the set
\begin{eqnarray*}\label{MG3}
\Lambda^I_\mu:=\{\lambda\in\mathbb{R}^+\dvtx \exists\sigma\in
\mathcal{C}^0(I)   \mbox{ with $(I,\mu,\sigma)$ is proper with FC
$\lambda$}\}
\end{eqnarray*}
and, if there is a $\lambda\in\Lambda^I_\mu$, we investigate
\begin{eqnarray*}\label{MG4}
\mathcal{H}^I_{\mu,\lambda}=\{\sigma\in\mathcal{C}^0(I)\dvtx
{}\mbox{$(I,\mu,\sigma)$ is proper with FC $\lambda$}\}.
\end{eqnarray*}
Proposition~5.1 of~\cite{KL} implies that there exist $I,\sigma$ with
$(I,\mu,\sigma)$ being proper if and only if $\beta<0$.
In this case it also follows that $I=\mathbb{R}$ and $\Lambda^I_\mu
=(0,|\beta|]$.
A FC $\lambda=|\beta|$ leads again to an affine model, namely $\sigma
(x)=\sigma_0 x$ for some $\sigma_0>0$.

If we choose an FC $\lambda=(1-\delta)|\beta|\in(0,|\beta|)$ for
some $\delta\in(0,1)$, then by Proposition~5.3 of~\cite{KL}, every
$\sigma\in\mathcal{H}^\mathbb{R}_{\mu,(1-\delta)|\beta|}$ is of
the form
\begin{eqnarray*}
\sigma(x)=\sigma_1|\alpha+\beta x|^\delta1_{\{x\leq
-\alpha/\beta\}}+\sigma_2|\alpha+\beta x|^\delta1_{\{
x\geq-\alpha/\beta\}}
\end{eqnarray*}
 for some $\sigma_1,\sigma_2>0$. Setting $f_i:=|\beta
|^{{\delta/(1-\delta)}}\sigma_i^{{1/(1-\delta)}}(1-\delta)^{
{1/(1-\delta)}}$ for $i=1,2$ the SST takes the form
\begin{eqnarray*}
f(x)=\biggl(\frac{\alpha}{\beta}-f_1|x|^{
{1/(1-\delta)}}\biggr)1_{\{x\leq0\}}+\biggl(\frac{\alpha}{\beta
}+f_2|x|^{{1/(1-\delta)}}\biggr)1_{\{x\geq0\}}.
\end{eqnarray*}
Calculating the derivative of $f,$ we see that a possible proper triple
is strongly proper if and only if $\delta\in[\frac{1}{2},1)$.
An example of such an SDE is for parameters ${\alpha}\in\mathbb{R}$
and $\beta
<0$ given by
%
\begin{eqnarray}\label{affineSDE}
\mathrm{d}X_t=(\alpha+\beta X_t)\,\mathrm{d}t+\sigma\sqrt{|\alpha+\beta X_t|}\,\mathrm{d}L^d_t
,\qquad t\in\mathbb{R}.
\end{eqnarray}
\end{example}

\begin{example}
Consider the SDE
%
\begin{eqnarray}\label{bss}
\mathrm{d}X_t=-\sigma_1\sin(\sigma_2 X_t)\cos(\sigma_2 X_t)\,\mathrm{d}t - \sin
^2(\sigma_2 X_t)\,\mathrm{d}L^d_t ,\qquad t\in\mathbb{R}.
\end{eqnarray}
This example provides a bounded state space model.
Define the triple $(I,\mu,\sigma)$ by $I:=(0,\frac{\curpi}{\sigma
_2})$, $\mu(x):=-\sigma_1\sin(\sigma_2 x)\cos(\sigma_2 x)$ and
$\sigma(x):=-\sin^2(\sigma_2 x)$ where $\sigma_1,\sigma_2>0$. It
can be shown that this triple is in fact strongly proper with FC
$\lambda=\sigma_1\sigma_2$. More precisely, we have $\psi(x)=\sigma
_1\cot(\sigma_2 x)$ and, therefore,
\begin{eqnarray*} X_t:=\frac{1}{\sigma_1}
\operatorname{arccot}(-\sigma_2 \mathcal{L}^{d,\sigma_1\sigma
_2}_t) ,\qquad t\in\mathbb{R},
\end{eqnarray*}
is the unique stationary solution of the SDE (\ref{bss}).
\end{example}

\subsection{Fractional Cox--Ingersoll--Ross models}\label{ss52}

Whenever positive phenomena are modeled -- for instance,
interest rates, volatilities or default rates in finance --
the Cox--Ingersoll--Ross (CIR)~\cite{coxinger} model is the most
prominent model.
It is the solution to
\begin{eqnarray*}
\mathrm{d}X_t= (a -\gamma X_t) \,\mathrm{d}t+\sigma\sqrt{X_t}\,\mathrm{d}B_t ,\qquad X_0=x_0\geq0,
\end{eqnarray*}
where $B=(B_t)_{t\in[0,\infty)}$ denotes standard Brownian motion,
$a,\gamma\in\mathbb{R}$ and ${\sigma}>0$.
General existence and uniqueness theorems of Brownian SDEs cannot be
applied here, because the square root is clearly not Lipschitz continuous.
However, Ikeda and Watanabe~\cite{IKWA}, page~221, showed that for any
$X_0=x\geq0$ there exists a unique non-negative solution.
We shall consider analogous SDEs driven by FLPs.

\begin{figure}[t]

\includegraphics{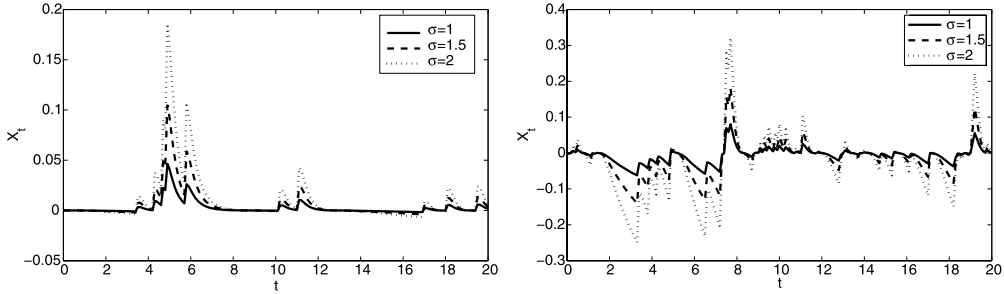}

\caption{Sample paths of a solution of the Cox--Ingersoll--Ross model
\protect\eqref{eq5.3} with $X_0=0$ for varying $\sigma$, fixed $\lambda=2.5$
and $d=0.35$, using two different FLP sample paths, left $\theta=0.5$,
right $\theta=2.5$.}\label{fig1}
\end{figure}

Within the framework of strongly proper triples,
Examples~\ref{exe1} and~\ref{exe2} show that our theory only covers
CIR models with mean reversion to $a=0$.
Consider for $\sigma,\gamma>0$ a solution to the pathwise SDE
%
\begin{equation}\label{eq5.3}
\mathrm{d}X_t=-\gamma X_t\,\mathrm{d}t+\sigma\sqrt{|X_t|}\,\mathrm{d}L^d_t ,\qquad t\in\mathbb{R}.
\end{equation}
Define $\widetilde{\sigma}(x):=\sigma|x|^{{1/2}}$,
choose $\widetilde{\mu}(x)=-\gamma x$ and take $I=\mathbb{R}$.
Example~\ref{exe1} implies that $(I,\widetilde{\mu},\widetilde
{\sigma})$ is strongly proper with SST
\begin{eqnarray*}
f(x)=\operatorname{sign}(x)\frac{{\sigma}^2}{4} x^2,
\end{eqnarray*}
and, by Theorem~\ref{th4.6}, a stationary solution of $(\ref{eq5.3})$
is given by $(f(\mathcal{L}^{d,\lambda}_t))_{t\in\mathbb{R}}$ with
$\lambda=\gamma/2$, cf.~Figure~\ref{fig1}.
Obviously, this CIR model takes also negative values.

%

\begin{figure}[t]

\includegraphics{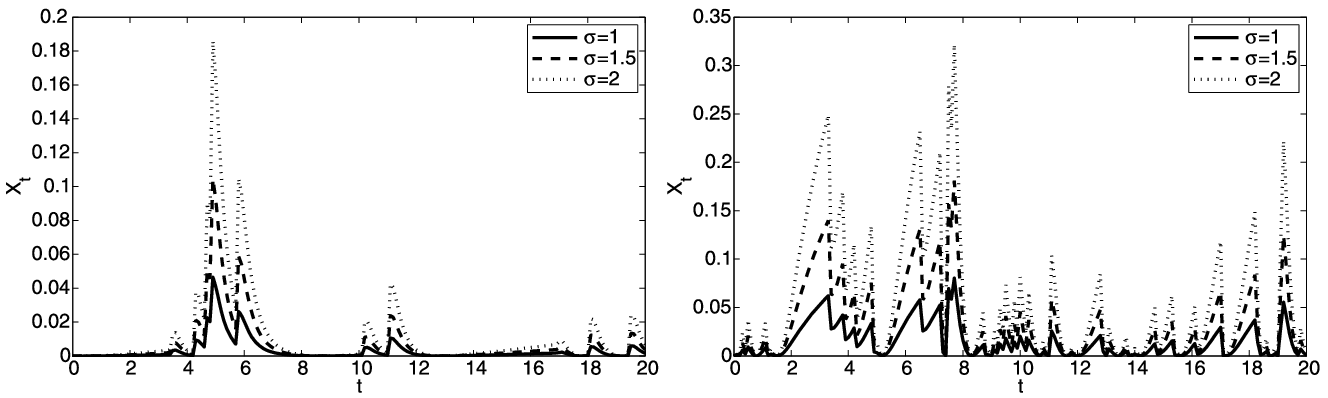}

\caption{Sample paths of squared FLOUPs for varying $\sigma$, fixed
$\lambda=2.5$ and $d=0.35$, using the same sample paths as in
Figure~\protect\ref{fig1}: left $\theta=0.5$, right $\theta=2.5$.}
\label{fig2}
\end{figure}

A natural non-negative transformation of the FLOUP is given by
$Z_t:=(\sigma\mathcal{L}_t^{d,\lambda})^2$ for $t\in\mathbb{R}$ (cf.
Figure~\ref{fig2}),
and, using the chain rule from Theorem~\ref{thA.3} and the existence
of all appearing Riemann--Stieltjes integrals, we get
\begin{eqnarray*}
\mathrm{d}Z_t&=&2\sigma^2\mathcal{L}^{d,\lambda}_t\,\mathrm{d}\mathcal{L}^{d,\lambda
}_t=-2\lambda\sigma^2(\mathcal{L}^{d,\lambda}_t)^2\,\mathrm{d}t+2\sigma
^2\mathcal{L}^{d,\lambda}_t\,\mathrm{d}L^d_t
\\
&=&-2\lambda Z_t\,\mathrm{d}t+2\sigma\sqrt
{Z_t}\,\mathrm{d}L^d_t,\qquad t\in\mathbb{R}.
\end{eqnarray*}
Defining now $\kappa(z):=-2\lambda z$ and $\iota(z):=2\sigma\sqrt
{z}$ we have $Z$ as a solution to
\begin{eqnarray*}\label{hui2}
\mathrm{d}Z_t=\kappa(Z_t)\,\mathrm{d}t+\iota(Z_t)\,\mathrm{d}L^d_t ,\qquad t\in\mathbb{R}.
\end{eqnarray*}
However, the triple $((0,\infty),\kappa,\iota)$ is not strongly
proper, because assumption (P2) of Definition~\ref{def4.1} is violated.

We can now formulate the following general result.

\begin{proposition}\label{pr5.4}
Let $\mathfrak{L}^{{\lambda/2}}(L^d,\cdot,\cdot)$ be the
Ornstein--Uhlenbeck operator from Definition~\ref{def3.9}.
Then for $\tau\in\mathbb{R}$ and $z\geq0$ the process
\begin{eqnarray*}
X_t^{\lambda,\tau,z}:=\biggl(\frac{\sigma}{2}\mathfrak{L}_t^{
{\lambda/2}}(L^d,\tau,z)\biggr)^2 ,\qquad t\in\mathbb{R},
\end{eqnarray*}
solves the SDEs
%
\begin{eqnarray}\label{eq5.13}
\mathrm{d}X_t&=&-\lambda X_t\,\mathrm{d}t+\sigma\sqrt{|X_t|}\,\mathrm{d}L^d_t
\end{eqnarray}
{and}
\begin{eqnarray}\label{eq5.14}
\mathrm{d}X_t&=&-\lambda X_t\,\mathrm{d}t+\sigma\sqrt{X_t}\,\mathrm{d}L^d_t ,\qquad t\in\mathbb
{R}.
\end{eqnarray}
In fact, any solution to \textup{(\ref{eq5.14})} also solves \textup{(\ref{eq5.13})}.
\end{proposition}

This result is not surprising, because Theorem~\ref{th4.7} does not
hold for the SDE \eqref{eq5.13}.
However, the reverse is not true: a solution of \eqref{eq5.13} does
not necessarily solve \eqref{eq5.14}, because it can be negative.
Also the constant process, $X_t:=0$, $t\in\mathbb{R}$, solves both
(\ref{eq5.14}) and (\ref{eq5.13}).

%

Considering a squared Ornstein--Uhlenbeck process leads in the case of
a driving Brownian motion to a CIR model with mean reversion to
positive values. This approach does not work for pathwise integrals,
neither in the FLP case nor for FBM, since the It\^o term in the chain
rule vanishes by finite $p$-variation for some $p<2$.

\begin{figure}[t]

\includegraphics{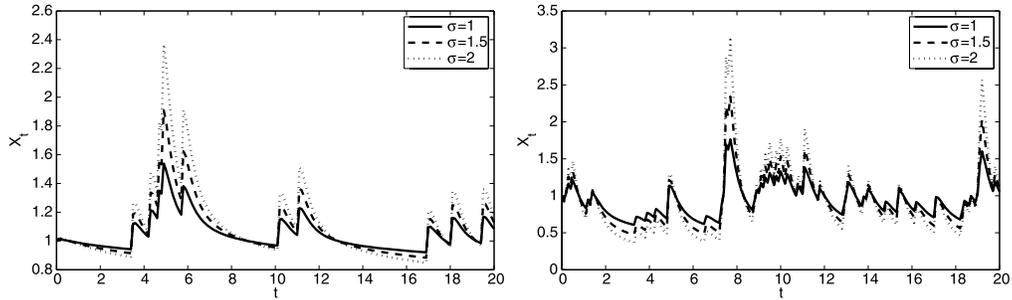}

\caption{Sample paths of a solution of (\protect\ref{eq5.6}) for varying
$\sigma$, fixed $\lambda=2.5$ and $d=0.35$, using the same sample
paths as in Figure~1: left $\theta=0.5$, right $\theta=2.5$.}
\label{fig3}
\end{figure}

A positive process based on Theorem~\ref{th4.6} is given as a solution to
%
\begin{equation}\label{eq5.6}
\mathrm{d}Y_t=-\lambda\sqrt{Y_t}\log(Y_t)\,\mathrm{d}t+\sigma
|Y_t|\,\mathrm{d}L^d_t,\qquad t\in\mathbb{R},\qquad Y_0=y_0
\end{equation}
for $\lambda,\sigma>0$, cf. Figure~\ref{fig3}.
Example~\ref{exe1} states that the triple for state space $I=(0,\infty
)$ is strongly proper and the SST can be calculated as $f(x)=\mathrm{e}^{\sigma x}$.

%

\begin{appendix}
\section*{Appendix: Riemann--Stieltjes integration}\label{app}

\renewcommand{\theequation}{A.\arabic{equation}}
\setcounter{equation}{0}
\renewcommand{\thetheorem}{A.\arabic{theorem}}
\setcounter{theorem}{0}

As mentioned in the \hyperref[s1]{Introduction}, all integrals in this paper are
considered Riemann--Stieltjes integrals, if not stated otherwise. That
is, for functions $f,h\dvtx [a,b]\mapsto\mathbb{R}$, we take the limit of
%
\begin{equation}\label{RSsums}
S(f,g,\kappa,\rho):=\sum_{i=1}^nf(y_i)[h(x_i)-h(x_{i-1})],
\end{equation}
where $\kappa=(x_i)_{i=0,\ldots,n}$ is a partition and $\rho
=(y_i)_{i=1,\ldots,n}$ an intermediate partition of $[a,b]$, that is,
\begin{eqnarray*}
a=x_0<x_1<\cdots<x_{n-1}<x_n=b,\qquad x_{i-1}\leq y_i\leq x_i\qquad
\mbox{for all }  i\in\{1,\ldots,n\},
\end{eqnarray*}
while letting $\operatorname{\textbf{mesh}}(\kappa):=\sup_{i=1,\ldots
,n}|x_i-x_{i-1}|$ go to zero. Using the Banach--Steinhaus theorem, one
can prove that if for a right-continuous $h$ and all continuous $f$ the
Riemann--Stieltjes sums of (\ref{RSsums}) converge, $h$ is already of
bounded variation. However, we can weaken this assumption on the
integrator by restricting the space of possible integrands. Recall the
definitions in \eqref{rev:1} and \eqref{wcon}. Exploiting the concept
of $p$-variation we now state an existence theorem for
Riemann--Stieltjes integrals proven by Young \cite{Young}.

\begin{theorem}\label{thA.2}
Let $[a,b]\subset\mathbb{R}$ be a compact interval, $f\in\mathfrak
{W}^{\mathrm{con}}_q([a,b])$ and $h\in\mathfrak{W}^{\mathrm{con}}_p([a,b])$
for some $p,q>0$ with $p^{-1}+q^{-1}>1$. Then $\int_a^bf_s\,\mathrm{d}h_s$ exists
in the Riemann--Stieltjes sense.
\end{theorem}

As in the classical Riemann--Stieltjes calculus, a chain rule can be
proven; see  \cite{za}, Theorem~3.1.

\begin{theorem}[(Chain rule)]\label{thA.3}
Let $[a,b]$ be a compact interval and $g\in\mathfrak{W}^{\mathrm
{con}}_p([a,b])$ for some $p\in(0,2)$. Furthermore,
let $F\in\mathcal{C}^1(\mathbb{R})$ with $F'\in\operatorname
{Lip}(\mathbb{R})$. Then the Riemann--Stieltjes integral $\int
_a^b(F'\circ g)_s\,\mathrm{d}g_s$
exists and we have
%
\begin{equation}\label{RS6}
(F\circ g)(b)-(F\circ g)(a)=\int_a^b(F'\circ g)_s\,\mathrm{d}g_s.
\end{equation}
\end{theorem}

At last we state a density formula, which we have not found in the
literature; for a proof we refer to  \cite{fink}, Theorem 4.3.2.

\begin{theorem}[(Density formula)]\label{thA.4}
Let $[a,b]\subset\mathbb{R}$ be a compact interval, $f,h\in\mathfrak
{W}^{\mathrm{con}}_q([a,b])$ and
$g\in\mathfrak{W}^{\mathrm{con}}_p([a,b])$ for some $q>0$ and $p>1$
with $p^{-1}+q^{-1}>1$. For all $x\in[a,b]$ we define $\phi(x):=\int
_a^xh_s\,\mathrm{d}g_s$.
Then we have $\phi\in\mathfrak{W}^{\mathrm{con}}_p([a,b])$ and
%
\begin{equation}\label{DF1}
\int_a^bf_s\,\mathrm{d}\phi_s=\int_a^bf_sh\,\mathrm{d}g_s.
\end{equation}
\end{theorem}
\end{appendix}

\section*{Acknowledgement}
We thank Martina Z\"ahle for interesting discussions and useful
comments, which led to an improvement of our paper.

\printhistory

\end{document}